\documentclass[a4paper,12pt,reqno]{amsart}

\usepackage{amsfonts}
\usepackage{amsmath}
\usepackage{amssymb}
\usepackage{cite}
\usepackage{mathrsfs}
\usepackage{enumitem}
\usepackage[colorlinks]{hyperref}

\numberwithin{equation}{section}

\usepackage{graphicx}

\newtheorem{theorem}{Theorem}[section]
\newtheorem{defi}[theorem]{Definition}

\def\R2n{{\mathbb R}^{2n}}

\def\R2{{\mathbb R}^2}
\def\R2n{{\mathbb R}^{2n}}

\def\N0{{\mathbb N}_{0}}

\def\l2h{{\ell^2(\hbar\mathbb Z^n)}}

\begin{document}

\title[Semiclassical Fractional Klein-Gordon equation]
{Discrete Time-dependent wave equations II. Semiclassical Fractional Klein-Gordon equation}

\author[Aparajita Dasgupta]{Aparajita Dasgupta}
\address{
	Aparajita Dasgupta:
	\endgraf
	Department of Mathematics
	\endgraf
	Indian Institute of Technology, Delhi, Hauz Khas
	\endgraf
	New Delhi-110016 
	\endgraf
	India
	\endgraf
	{\it E-mail address} {\rm adasgupta@maths.iitd.ac.in}
}
\author[Michael Ruzhansky]{Michael Ruzhansky}
\address{
	Michael Ruzhansky:
	\endgraf
	Department of Mathematics: Analysis, Logic and Discrete Mathematics
	\endgraf
	Ghent University, Belgium
	\endgraf
	and
	\endgraf
	School of Mathematical Sciences
	\endgraf
	Queen Mary University of London
	\endgraf
	United Kingdom
	\endgraf
	{\it E-mail address} {\rm ruzhansky@gmail.com}
}
\author[Abhilash Tushir]{Abhilash Tushir}
\address{
	Abhilash Tushir:
	\endgraf
	Department of Mathematics
	\endgraf
	Indian Institute of Technology, Delhi, Hauz Khas
	\endgraf
	New Delhi-110016 
	\endgraf
	India
	\endgraf
	{\it E-mail address} {\rm abhilash2296@gmail.com}
}

\thanks{ The first author was supported by Core Research Grant, RP03890G,  Science and Engineering
	Research Board,  India. The second
 author was supported by the EPSRC Grant 
EP/R003025/2, by the FWO Odysseus 1 grant G.0H94.18N: Analysis and Partial Differential Equations and by the  Methusalem programme of the Ghent University Special Research Fund (BOF) (Grantnumber
01M01021).}
\date{\today}


\begin{abstract}
 In this paper we consider a semiclassical version of the fractional Klein-Gordon equation  on the lattice $\hbar \mathbb{Z}^{n} .$ Contrary to the Euclidean case that was considered in \cite{fracklien}, the discrete fractional Klein-Gordon equation  is well-posed in $\ell^{2}\left(\hbar \mathbb{Z}^{n}\right) .$ However, we also recover the well-posedness results in the  certain
  Sobolev spaces in the limit of the semiclassical parameter $\hbar\to 0$. 
\end{abstract}
\maketitle

\tableofcontents
\section{Introduction}
In this paper we study the fractional Klein-Gordon equation on the discrete lattice
$$
\hbar \mathbb{Z}^{n}=\left\{x \in \mathbb{R}^{n}: x=\hbar k, k \in \mathbb{Z}^{n}\right\},
$$
depending on a (small) discretization parameter $\hbar>0$, and the behaviour of its solutions as $\hbar \rightarrow 0.$ The discrete fractional Laplacian on $\hbar \mathbb{Z}^{n}$ is denoted by $\left(-\mathcal{L}_{\hbar}\right)^{\alpha}$ and is defined by
\begin{equation}\label{fraclap}
	\left(-\mathcal{L}_{\hbar}\right)^{\alpha}u(k):=\sum\limits_{j\in\mathbb{Z}^{n}}a^{(\alpha)}_{j}u(k+j\hbar),\quad 0<\alpha<1,
\end{equation}
where the expansion coefficient $a_{j}^{(\alpha)}$ is given by
\begin{equation}
		a_{j}^{(\alpha)}:=\int\limits_{\left[-\frac{1}{2}
			,\frac{1}{2}\right]^{n}} \left|2\sin\left(\pi\theta\right)\right|^{2\alpha}e^{-2\pi ij\cdot\theta}\mathrm{d}\theta.
\end{equation}
 We refer to Section \ref{prlim} for further details about the discrete fractional Laplacian.

The fractional Klein-Gordon equation with fractional Laplacian on $\mathbb{R}^{n}$ is given by
\begin{equation}\label{orgpde}
	\left\{\begin{array}{l}
		\partial_{t}^{2} u(t, x)+ (-\mathcal{L})^{\alpha} u(t, x)+m(x) u(t, x)=f(t,x), \quad \text{ with }  (t,x)\in(0,T]\times\mathbb{R}^{n},\ \\
		u(0, x)=u_{0}(x),\quad x\in\mathbb{R}^{n}, \\
		\partial_{t} u(0, x)=u_{1}(x),\quad x\in\mathbb{R}^{n},
	\end{array}\right.
\end{equation}
where $\left(-\mathcal{L}\right)^{\alpha}$ is the usual fractional Laplacian on $\mathbb{R}^{n}$  defined as a pseudo-differential operator with symbol $|2\pi\xi|^{2\alpha}=\left[\sum\limits_{j=1}^{n}4\pi^{2}\xi_{j}^{2}\right]^{\alpha}$ i.e.,
\begin{equation}\label{frlap}
	(-\mathcal{L})^{\alpha} u(x)=\mathcal{F}^{-1}\left\{\left|2\pi\xi\right|^{2\alpha} \mathcal{F} u(\xi)\right\}\left(x\right)=\int\limits_{\mathbb{R}^{n}}e^{2\pi ix\cdot\xi}\left|2\pi\xi\right|^{2\alpha}\widehat{u}(\xi)\mathrm{d}\xi,
\end{equation}
where $\mathcal{F}$ and $\mathcal{F}^{-1}$ are the usual Fourier transform and inverse Fourier transform, respectively. The fractional Laplacian operator can also be understood as a positive power of the classical Laplacian. We refer to \cite{fracthought} and \cite{ten} for more details and other alternative definitions of the fractional Laplacian.

The nonlocal nature of fractional derivatives like  Riemann–Liouville fractional derivative, Caputo fractional derivative, Riesz fractional derivative, etc., plays a vital role in modeling various problems of  classical and quantum mechanics  (see \cite{gol,gepr}). In the last few decades, considerable attention has been given to the solutions of fractional partial differential equations. 

Recently, analytical and numerical solutions of the fractional Klein-Gordon equation have been studied by many authors like \cite{abuteen,culha,ege,gepr,garra}.
  There have been studies where the mass term depends on the position, \cite{altu,desouza,wang,wang1} and in \cite{ghosh,gol} the fractional Laplacian was introduced. Recently in \cite{fracklien}, Altybay et al. considered the case where the mass term is singular with zero source term and studied the well-posedness in the very weak sense.  In the regular situation i.e., when the mass term $m$ is a regular function, authors have proved that the  Cauchy problem \eqref{orgpde} has a unique solution $u \in C\left([0, T] ; H^{\alpha}\left(\mathbb{R}^{n}\right)\right) \cap C^{1}\left([0, T] ; L^{2}\left(\mathbb{R}^{n}\right)\right)$, where $H^{\alpha}(\mathbb{R}^{n})$ is a fractional Sobolev space. More precisely, they have obtained the following well-posedness result for $f\equiv 0$ : 
\medskip
\begin{theorem}\label{fracklien}
	Let  $ m \in L^{\infty}\left(\mathbb{R}^{n}\right)$ be non-negative. Suppose that $u_{0} \in H^{\alpha}\left(\mathbb{R}^{n}\right)$ and $u_{1} \in  L^{2}\left(\mathbb{R}^{n}\right)$. Then, there is a unique solution $u \in C\left([0, T] ; H^{\alpha}\left(\mathbb{R}^{n}\right)\right) \cap C^{1}\left([0, T] ; L^{2}\left(\mathbb{R}^{n}\right)\right)$ to \eqref{orgpde}, and it satisfies the estimate
	\begin{equation}
		\left\|u(t,\cdot)\right\|^{2}_{H^{\alpha}(\mathbb{R}^{n})}+\left\|\partial_{t}u(t,\cdot)\right\|^{2}_{L^{2}(\mathbb{R}^{n})} \lesssim\left(1+\|m\|_{L^{\infty}(\mathbb{R}^{n})}\right)\left[\left\|u_{0}\right\|_{H^{\alpha}(\mathbb{R}^{n})}^{2}+\left\|u_{1}\right\|_{L^{2}(\mathbb{R}^{n})}^{2}\right],
	\end{equation}
	where $	\left\|u\right\|^{2}_{H^{\alpha}(\mathbb{R}^{n})}=\left\|u\right\|^{2}_{L^{2}(\mathbb{R}^{n})}+\left\|(-\mathcal{L})^{\frac{\alpha}{2}}u\right\|^{2}_{L^{2}(\mathbb{R}^{n})}$.
\end{theorem}

Furthermore in \cite{mari}, Chatzakou et al. also studied the Klein-Gordon  equation with positive (left) Rockland operator  on a general graded Lie group. This setting of Rockland operators on graded Lie groups allows one to consider both elliptic and sub-elliptic settings in \ref{orgpde}. In this paper we are interested in studying the discretization of the Klein-Gordon equation with the fractional Laplacian like the classical wave equation in \cite{dasgupta2021discrete}. \\

The semiclassical analogue of the fractional Klein-Gordon equation  on the lattice $\hbar \mathbb{Z}^{n}$  is given by the Cauchy problem
\begin{equation}\label{mainpde}
	\left\{\begin{array}{l}
		\partial_{t}^{2} u(t, k)+\hbar^{-2\alpha}(-\mathcal{L}_{\hbar})^{\alpha} u(t, k)+m(k) u(t, k)=f(t,k), \quad \text { with } (t,k) \in(0, T]\times\hbar\mathbb{Z}^{n}, \\
		u(0, k)=u_{0}(k), \quad k \in \hbar \mathbb{Z}^{n}, \\
		\partial_{t} u(0, k)=u_{1}(k), \quad k \in \hbar \mathbb{Z}^{n},
	\end{array}\right.
\end{equation}
 with non-negative mass term  $m\in\ell^{\infty}(\hbar\mathbb{Z}^{n})$ and  $f\in L^{2}([0,T];\ell^{2}(\hbar\mathbb{Z}^{n}))$. 
The well-posedness of the  Cauchy problem  \eqref{mainpde} on the discrete lattice $\hbar\mathbb{Z}^{n}$ varies from the  results obtained in Theorem \ref{fracklien} for the Euclidean case
in the sense that \eqref{mainpde} is always well-posed in $\ell^{2}(\hbar\mathbb{Z}^{n})$, while in the Euclidean case, it is well-posed in fractional Sobolev space. More precisely, we shall prove the following well-posedness result for the lattice case in Section \ref{secclass}.
\medskip
\begin{theorem}\label{wellpo}
	Let $T>0$. Assume that $m \in \ell^{\infty}\left(\hbar\mathbb{Z}^{n}\right)$  and $f\in L^{2}([0,T];\ell^{2}(\hbar\mathbb{Z}^{n}))$, then the Cauchy problem \eqref{mainpde} is well-posed in $\ell^{2}(\hbar\mathbb{Z}^{n})$. In particular, if $u_{0},u_{1} \in \ell^{2}\left(\hbar\mathbb{\mathbb{Z}}^{n}\right)$ and $f\in L^{2}([0,T];\ell^{2}(\hbar\mathbb{Z}^{n}))$, then for every $t\in[0,T]$, we have $u(t),\partial_{t}u(t) \in  \ell^{2}\left(\hbar\mathbb{Z}^{n}\right)$. Moreover for each $\hbar>0$, it satisfies the estimate
	\begin{multline}\label{uestt}
	\|u(t, \cdot)\|^{2}_{\ell^{2}(\hbar\mathbb{Z}^{n})}+\left\|\partial_{t} u(t, \cdot)\right\|_{\ell^{2}(\hbar\mathbb{Z}^{n})}^{2}\lesssim\left(1+\|m\|_{\ell^{\infty}(\hbar\mathbb{Z}^{n})}\right)\times\\\left[\left(\hbar^{-2\alpha}+\|m\|_{\ell^{\infty}(\hbar\mathbb{Z}^{n})}\right)\left\|u_{0}\right\|_{\ell^{2}(\hbar\mathbb{Z}^{n})}^{2}+\left\|u_{1}\right\|_{\ell^{2}(\hbar\mathbb{Z}^{n})}^{2}+\|f\|^{2}_{L^{2}([0,T];\ell^{2}(\hbar\mathbb{Z}^{n}))}\right],
	\end{multline}	
 for all $u_{0},u_{1}\in\ell^{2}(\hbar\mathbb{Z}^{n})$ and $f\in L^{2}([0,T];\ell^{2}(\hbar\mathbb{Z}^{n}))$. The constant depends on $T$ but not on  $\hbar$.
\end{theorem}
\medskip
A natural question one may ask here is:
\begin{itemize}
	\item[] Why is there a difference in the well-posedness results between the Euclidean case and the lattice case?
\end{itemize}
In view of  Theorem \ref{wellpo},  the dependence of RHS of estimate \eqref{uestt} on $\hbar$ can answer this question immediately. Clearly, RHS may go to infinity as $\hbar \rightarrow 0$. This difference arises because of the difference in boundedness behavior of the discrete fractional  Laplacian  and the fractional Laplacian on  $\ell^{2}(\hbar\mathbb{Z}^{n})$ and $L^{2}(\mathbb{R}^{n})$, respectively. 

Further, we are also interested in  approximating the continuous solution in Euclidean settings
 by the discrete solution in the lattice settings. There are many ways to approximate but a large literature includes the numerical approaches or by minimizing the notion of distance between the continuous solution and the discrete one, see for example \cite{hua,kirk,roncal} and references therein. The nonlocality and singularity of the fractional Laplace operator in equation \eqref{orgpde} are major obstacles in numerical approach. So, in  this paper we are interested in approximating the continuous solution by the discrete solution in  $\ell^{2}$-norm  like  we did for the classical wave equation in \cite{dasgupta2021discrete}. The following theorem shows that under the assumptions that the solutions in the Euclidean case exist, they can be approximated/recovered in the semiclassical  limit as $\hbar \rightarrow 0$. We require a little additional Sobolev
regularity to ensure that the convergence results are global on the whole of $\mathbb{R}^{n}$. 
\begin{theorem}\label{cgt}
Let $u$ and $v$ be the solutions of the Cauchy problems \eqref{mainpde} on $\hbar \mathbb{Z}^{n}$ and \eqref{orgpde} on $\mathbb{R}^{n}$, respectively. Assume that the Cauchy data $u_{0}\in H^{4\alpha}(\mathbb{R}^{n})$ and $u_{1}\in L^{2}(\mathbb{R}^{n})$, then for every $t \in[0, T]$, we have
\begin{equation}
\|v(t)-u(t)\|_{\ell^{2}\left(\hbar \mathbb{Z}^{n}\right)}+\left\|\partial_{t} v(t)-\partial_{t} u(t)\right\|_{\ell^{2}\left(\hbar \mathbb{Z}^{n}\right)} \rightarrow 0\text{ as } \hbar \rightarrow 0 .
\end{equation}

\end{theorem}
Thus, when $\hbar\to 0$, for the above mentioned space, we actually recover the well-posedness
results on $\mathbb{R}^{n}$, where the solution $u$ in the above statement becomes restricted to the lattice $\hbar\mathbb{Z}^{n}$.

We note that the symbolic calculus of pseudo-difference operators on lattice $\mathbb{Z}^{n}$ has been developed  in \cite{kibti}. The symbolic calculus on $\mathbb{Z}^{n}$ can be thought
of as a dual one to the calculus developed on the torus in \cite{RT-Birk,RT-JFAA}. The difference equations on lattice $\hbar\mathbb{Z}^{n}$, including Schrödinger equations have been
studied (e.g. in \cite{Rab13,Rab10,markus,Rab18}) by developing the analysis in terms of kernels.
\begin{defi}[Symbol class $S^{m}\left(\hbar \mathbb{Z}^{n} \times \mathbb{T}_{\hbar}^{n}\right)$]
	Let $m \in(-\infty, \infty)$. We say that a function $\sigma$ : $\hbar \mathbb{Z}^{n} \times \mathbb{T}_{\hbar}^{n} \rightarrow \mathbb{C}$ belongs to $S^{m}\left(\hbar \mathbb{Z}^{n} \times \mathbb{T}_{\hbar}^{n}\right)$ if $\sigma(k, \cdot) \in C^{\infty}\left(\mathbb{T}_{\hbar}^{n}\right)$ for all $k \in \hbar \mathbb{Z}^{n}$, and for all multi-indices $\alpha, \beta$ there exists a positive constant $C_{\alpha, \beta, \hbar}$ such that
	\begin{equation}
		\left|D_{\theta}^{(\beta)} \Delta_{k}^{\alpha} \sigma(k, \theta)\right| \leq C_{\alpha, \beta, \hbar}(1+|k|)^{m-|\alpha|},
	\end{equation} 
	for all $k \in \hbar \mathbb{Z}^{n}$ and $\theta \in \mathbb{T}_{\hbar}^{n}$. We denote by $\mathrm{Op}(\sigma)$ the pseudo-difference operator with symbol $\sigma$ given by
	\begin{equation}
		\operatorname{Op}(\sigma) f(k):=\hbar^{n/2}\int_{\mathbb{T}_{\hbar}^{n}} e^{2 \pi i k \cdot \theta} \sigma(k, \theta) \widehat{f}(\theta) \mathrm{d} \theta, \quad k \in \hbar \mathbb{Z}^{n}
	\end{equation}
	where
	\begin{equation}
		\widehat{f}(\theta)=\hbar^{n/2}\sum_{k \in \hbar \mathbb{Z}^{n}}f(k) e^{-2 \pi i k \cdot \theta} , \quad \theta \in \mathbb{T}_{\hbar}^{n}.
	\end{equation}
\end{defi}

This paper is arranged in the following four sections. In Section \ref{prlim}, we develop basic Fourier analysis  and   discrete fractional Laplacian on lattice $\hbar\mathbb{Z}^{n}$. 
In Section \ref{secclass}, we establish the  well-posedness result for the Cauchy problem \eqref{mainpde}. Finally, in
Section \ref{limit},  we discuss the limiting behaviour of
solutions to \eqref{mainpde} in the limit of the semiclassical parameter $\hbar\to 0$.

To simplify the notation, throughout the paper we will be writing $A\lesssim B$ if there
exists a constant $C$ independent of the appearing parameters such that $A\leq CB$.
\section{Preliminaries}\label{prlim}
 In this section, we first recall some facts concerning the Fourier analysis on the discrete lattice $\hbar\mathbb{Z}^{n}$.

The Schwartz space $\mathcal{S}\left(\hbar \mathbb{Z}^{n}\right)$ on the lattice $\hbar \mathbb{Z}^{n}$ is the space of rapidly decreasing functions $u: \hbar \mathbb{Z}^{n} \rightarrow \mathbb{C}$, that is, $u \in \mathcal{S}\left(\hbar \mathbb{Z}^{n}\right)$ if for any $L<\infty$ there exists a constant $C_{u, L, \hbar}$ such that
\begin{equation}
|u(k)| \leq C_{u, L, \hbar}(1+|k|)^{-L}, \quad \text { for all } k \in \hbar \mathbb{Z}^{n},
\end{equation}
where $|k|=\hbar\left(\sum\limits_{j=1}^{n} k_{j}^{2}\right)^{1 / 2} .$ A grid function $u: \hbar\mathbb{Z}^{n} \rightarrow \mathbb{R}$ is in $\ell^{p}(\hbar\mathbb{Z}^{n}),$ $ 1 \leq p<\infty$,  if
\begin{equation}
	\|u\|_{\ell^{p}(\hbar\mathbb{Z}^{n})}:=\left( \sum_{k \in \hbar\mathbb{Z}^{n}}\left|u(k)\right|^{p}\right)^{1 / p}<\infty,
\end{equation}
while $u \in \ell^{\infty}(\hbar\mathbb{Z}^{n})$ if
\begin{equation}
\|u\|_{\ell^{\infty}(\hbar\mathbb{Z}^{n})}:=\sup _{k \in \hbar\mathbb{Z}^{n}}\left|u(k)\right|<\infty .
\end{equation}
If $1 \leq p < q \leq \infty$, then 
\begin{equation}
\ell^{p}(\hbar\mathbb{Z}^{n}) \subsetneq \ell^{q}(\hbar\mathbb{Z}^{n})\text{ and }\|u\|_{\ell^{q}(\hbar\mathbb{Z}^{n})} \lesssim \|u\|_{\ell^{p}(\hbar\mathbb{Z}^{n})}.
\end{equation}
 Let $u\in\ell^{p}(\hbar\mathbb{Z}^{n})$  and $v\in\ell^{q}(\hbar\mathbb{Z}^{n})$, with $1\leq p,q\leq \infty$ satisfying  $\frac{1}{p}+\frac{1}{q}=1$, then the Hölder's inequality takes the form
\begin{equation}
\|uv\|_{\ell^{1}(\hbar\mathbb{Z}^{n})} \leq\|u\|_{\ell^{p}(\hbar\mathbb{Z}^{n})}\|v\|_{\ell^{q}(\hbar\mathbb{Z}^{n})}.
\end{equation}
The inner product in the Hilbert space $\ell^{2}\left(\hbar \mathbb{Z}^{n}\right)$ is given by
\begin{equation}
	\left( u,v\right)_{\ell^{2}\left(\hbar \mathbb{Z}^{n}\right)}=\sum_{k \in \hbar \mathbb{Z}^{n}} u(k) \overline{v(k)},
\end{equation}
where $\overline{v}$ is the complex conjugate of $v$. Let $\mathbb{T}_{\hbar}^{n}$ be the  $n$-dimensional torus, which we identify with periodic  $\left[-\frac{1}{2\hbar},\frac{1}{2\hbar}\right]^{n}$.
The Fourier transform $\mathcal{F}_{\hbar\mathbb{Z}^{n}}:\ell^{2}(\hbar\mathbb{Z}^{n})\to L^{2}(\mathbb{T}_{\hbar}^{n})$  is defined as
\begin{equation}
\mathcal{F}_{\hbar\mathbb{Z}^{n}}u(\theta):=\hbar^{n/2}\sum_{k \in \hbar \mathbb{Z}^{n}} u(k) e^{-2 \pi i k \cdot \theta}, \quad u\in \ell^{2}(\hbar\mathbb{Z}^{n}) ,
\end{equation}
and the inverse Fourier transform $\mathcal{F}^{-1}_{\hbar\mathbb{Z}^{n}}:L^{2}(\mathbb{T}_{\hbar}^{n})\to\ell^{2}(\hbar\mathbb{Z}^{n})$ is defined as
\begin{equation}
	\mathcal{F}^{-1}_{\hbar\mathbb{Z}^{n}}v(k):=\hbar^{n/2}\int_{\mathbb{T}_{\hbar}^{n}}v(\theta)e^{2\pi i k\cdot\theta}\mathrm{d}\theta, \quad v\in L^{2}(\mathbb{T}_{\hbar}^{n}).
\end{equation}
The Plancherel formula takes the form
 \begin{eqnarray}
\int_{\mathbb{T}_{\hbar}^{n}}|\widehat{u}(\theta)|^{2} d \theta&=&\int_{\mathbb{T}_{\hbar}^{n}} \widehat{u}(\theta) \overline{\widehat{u}(\theta)} d \theta\nonumber\\
 		&=&\int_{\mathbb{T}_{\hbar}^{n}} \hbar^{n}\sum_{k \in \hbar \mathbb{Z}^{n}} u(k) e^{-2\pi i k \cdot \theta} \sum_{m \in \hbar \mathbb{Z}^{n}} \overline{u(m)} e^{2\pi i m \cdot \theta} d \theta\nonumber\\
 		&=&\hbar^{n}\sum_{k \in \hbar \mathbb{Z}^{n}}\sum_{m \in \hbar \mathbb{Z}^{n}}u(k)\overline{u(m)}\int_{\mathbb{T}_{\hbar}^{n}}e^{2\pi i(m-k)\cdot\theta}\mathrm{d}\theta\nonumber\\
 		&=&\sum_{k \in \hbar \mathbb{Z}^{n}}|u(k)|^{2},
 	 \end{eqnarray}
  and the  Parseval's identity 
\begin{equation}
	\left( u,v\right)_{\ell^{2}\left(\hbar \mathbb{Z}^{n}\right)}=	\int_{\mathbb{T}_{\hbar}^{n}}\widehat{u}(\theta)\overline{\widehat{v}(\theta)} \mathrm{d} \theta.
\end{equation}
 The Fourier inversion formula is given by
\begin{equation}
u(k)=\hbar^{n/2}\int_{\mathbb{T}_{\hbar}^{n}}\widehat{u}(\theta) e^{2 \pi i k \cdot \theta}  \mathrm{d}\theta, \quad k \in \hbar\mathbb{Z}^{n}.
\end{equation}


Let $\alpha>0$ and let $u(x)$ be a complex-valued function. The fractional centered difference (FCD) operator is given by
\begin{equation}
	\Delta_{c}^{\alpha} u(x):=\sum\limits_{j\in\mathbb{Z}} \frac{(-1)^{j} \Gamma(\alpha+1)}{\Gamma(\alpha / 2-j+1) \Gamma(\alpha / 2+j+1)} u(x+j \hbar).
\end{equation}
Using the following relation from \cite[page 114]{duarte}:
\begin{equation}
	|2 \sin (\pi\theta)|^{2\alpha}=\sum\limits_{j\in\mathbb{Z}}  \frac{(-1)^{j} \Gamma(\alpha+1)}{\Gamma(\alpha / 2-j+1) \Gamma(\alpha / 2+j+1)} e^{2\pi i j\cdot\theta}, 
\end{equation}
we get
\begin{equation}
	 \frac{(-1)^{j} \Gamma(\alpha+1)}{\Gamma(\alpha / 2-j+1) \Gamma(\alpha / 2+j+1)}=\int_{-\frac{1}{2}}^{\frac{1}{2}} |2\sin(\pi\theta)|^{2\alpha}e^{-2\pi ij\cdot\theta}\mathrm{d}\theta,
\end{equation}
where $|2\sin(\pi\theta)|^{2\alpha}=\left[\sum\limits_{i=1}^{n}4\sin^{2}\left(\pi\theta_{i}\right)\right]^{\alpha}$ Since, the integral form of the above identity make sense for $j\in \mathbb{Z}^{n}$ as well, therefore, the definition of the fractional centered difference can be extended for functions on $\mathbb{R}^{n}$ as well in the following way
\begin{equation}\label{fracdef}
	\Delta_{c}^{\alpha} u(x):=\sum\limits_{j\in\mathbb{Z}^{n}} a_{j}^{(\alpha)} u(x+j \hbar),
\end{equation}
where the generating function $a_{j}^{(\alpha)}$ is given by
\begin{equation}\label{akalpha}
	a_{j}^{(\alpha)}:=\int\limits_{\left[-\frac{1}{2},\frac{1}{2}\right]^{n}} |2\sin(\pi\theta)|^{2\alpha}e^{-2\pi ij\cdot\theta}\mathrm{d}\theta=\int\limits_{\left[-\frac{1}{2}, \frac{1}{2}\right]^{n}}\left[\sum\limits_{i=1}^{n}4\sin^{2}\left(\pi\theta_{i}\right)\right]^{\alpha}e^{-2\pi ij\cdot\theta}\mathrm{d}\theta.
\end{equation}
Using the Fourier transform and the Fourier inversion formula, it is very easy to verify that
\begin{equation}
	|2 \sin (\pi\theta)|^{2\alpha}=\sum\limits_{j\in\mathbb{Z}^{n}} a_{j}^{(\alpha)} e^{2\pi i j\theta}.
\end{equation}
For more details about the FCD, we refer to \cite{duarte,samko1993fractional}. 

The discrete fractional Laplacian on $\hbar\mathbb{Z}^{n}$ is defined using the  FCD. The discrete fractional Laplacian on $\hbar \mathbb{Z}^{n}$,  denoted by $\left(-\mathcal{L}_{\hbar}\right)^{\alpha}$,  is defined by
\begin{equation}
	\left(-\mathcal{L}_{\hbar}\right)^{\alpha}u(k):=\sum\limits_{j\in\mathbb{Z}^{n}}a^{(\alpha)}_{j}u(k+j\hbar),\quad 0<\alpha<1,
\end{equation}
where $a^{(\alpha)}_{j}$ is  given by \eqref{akalpha}.
 For example,  take
$\alpha=1$, then it is very easy to check that
\begin{itemize}
	\item[1)] for $j=0$, we have $a_{j}^{(1)}=2n$,
	\item[2)] for $j=v_{i}$, we have $a_{j}^{(1)}=-1$,
	\item[3)] for $j=-v_{i}$, we have $a_{j}^{(1)}=-1$,
	\item[4)] for $j\neq 0,\pm v_{i}$, we have $a_{j}^{(1)}=0$,
\end{itemize}
where  $v_{i}$  is the  $i^{t h}$  basis vector in  $\mathbb{Z}^{n}$, having all zeros except for  $1$  at the  $i^{t h}$  component. 
This gives
\begin{eqnarray}
	(-\mathcal{L}_{\hbar})^{1}u(k)&=&2nu(k)-\sum\limits_{j=\pm v_{i}}u(k+j\hbar)\nonumber\\
	&=&2nu(k)-\sum\limits_{i=1}^{n}\left(u(k+v_{i}\hbar)+u(k-v_{i}\hbar)\right),
\end{eqnarray}
which is usual discrete Laplacian on $\hbar\mathbb{Z}^{n}$. For more details about the discrete Laplacian on $\hbar\mathbb{Z}^{n}$, one can refer to \cite{dasgupta2021discrete}.

 The symbol of $(-\mathcal{L}_{\hbar})^{\alpha}$ defined by  $\sigma_{(-\mathcal{L}_{\hbar})^{\alpha}}(k,\theta)=e^{-2\pi i k\cdot \theta}(-\mathcal{L}_{\hbar})^{\alpha}e^{2\pi i k\cdot \theta}$ is given by 
\begin{eqnarray}\label{fcdf}
	\sigma_{(-\mathcal{L}_{\hbar})^{\alpha}}(k,\theta)&=&e^{-2\pi i k\cdot \theta}(-\mathcal{L}_{\hbar})^{\alpha}e^{2\pi i k\cdot \theta}\nonumber\\
	&=&e^{-2\pi i k\cdot \theta}\sum\limits_{j\in\mathbb{Z}^{n}} a_{j}^{(\alpha)} e^{2\pi i (k+j \hbar)\cdot \theta}\nonumber\\
	&=&\sum\limits_{j\in\mathbb{Z}^{n}} a_{j}^{(\alpha)}e^{2\pi i j\hbar\cdot \theta}\nonumber\\
	&=&|2 \sin (\pi\hbar\theta)|^{2\alpha},
\end{eqnarray}
 with $(k, \theta) \in \hbar \mathbb{Z}^{n} \times \mathbb{T}_{\hbar}^{n}$, and it is independent of  $k$. 

  Moreover, since the symbol $\left|2\sin\left(\pi\hbar\theta\right)\right|^{2\alpha}$  is non-negative, therefore the discrete fractional Laplacian is a self-adjoint operator i.e.,
\begin{equation}
	((-\mathcal{L}_{\hbar})^{\alpha} u, u)_{\ell^{2}(\hbar\mathbb{Z}^{n})}=((-\mathcal{L}_{\hbar})^{\frac{\alpha}{2}} u,(-\mathcal{L}_{\hbar})^{\frac{\alpha}{2}} u)_{\ell^{2}(\hbar\mathbb{Z}^{n})}=( u,(-\mathcal{L}_{\hbar})^{\alpha} u)_{\ell^{2}(\hbar\mathbb{Z}^{n})}.
\end{equation}
\section{Proof of Theorem \ref{wellpo}}\label{secclass}
In this section, we will study the existence and uniqueness of the solution of the Cauchy problem \eqref{mainpde}.

\begin{proof}[Proof of Theorem \ref{wellpo}]
	Considering the discrete inner product of  \eqref{mainpde} with $\partial_{t}u$ on both sides and comparing the real part, we get
	\begin{multline}\label{wellm}
		\operatorname{Re}\left(\left(\partial_{t}^{2} u(t, \cdot), \partial_{t} u(t, \cdot)\right)_{\ell^{2}(\hbar\mathbb{Z}^{n})}+\left(\hbar^{-2\alpha}(-\mathcal{L}_{\hbar})^{\alpha} u(t, \cdot), \partial_{t} u(t, \cdot)\right)_{\ell^{2}(\hbar\mathbb{Z}^{n})} 
		+\right.\\
		\left.\left( m(\cdot) u(t, \cdot), \partial_{t} u(t, \cdot)\right)_{\ell^{2}(\hbar\mathbb{Z}^{n})}\right)=\operatorname{Re}\left(  f(t, \cdot), \partial_{t} u(t, \cdot)\right)_{\ell^{2}(\hbar\mathbb{Z}^{n})},
	\end{multline}
	for all $t\in(0,T]$. 
	Using  the self-adjointness of $(-\mathcal{L}_{\hbar})^{\alpha}$ and the techniques from \cite{fracklien,fracschro}, it is easy to check that:
	\begin{itemize}
		\item[i)] $\operatorname{Re}\left(\partial_{t}^{2} u(t, \cdot), \partial_{t} u(t, \cdot)\right)_{\ell^{2}(\hbar\mathbb{Z}^{n})}=\frac{1}{2}\partial_{t}\left\|\partial_{t} u(t,\cdot)\right\|^{2}_{\ell^{2}(\hbar\mathbb{Z}^{n})}$,\\
		\item[ii)] $\operatorname{Re}\left(\hbar^{-2\alpha}(-\mathcal{L}_{\hbar})^{\alpha} u(t, \cdot), \partial_{t} u(t, \cdot)\right)_{\ell^{2}(\hbar\mathbb{Z}^{n})}=\frac{1}{2}\partial_{t}\left\|\hbar^{-\alpha}(-\mathcal{L}_{\hbar})^{\frac{\alpha}{2}} u(t,\cdot)\right\|^{2}_{\ell^{2}(\hbar\mathbb{Z}^{n})}$,\\
		\item[iii)] $\operatorname{Re}\left( m(\cdot) u(t, \cdot), \partial_{t} u(t, \cdot)\right)_{\ell^{2}(\hbar\mathbb{Z}^{n})}=\frac{1}{2}\partial_{t}\left\|m^{\frac{1}{2}}(\cdot) u(t,\cdot)\right\|^{2}_{\ell^{2}(\hbar\mathbb{Z}^{n})}.$
	\end{itemize}
	
The energy functional of the system \eqref{mainpde} is defined by
	\begin{equation}\label{enrgy}
		E(t):=\left\|\partial_{t} u(t, \cdot)\right\|_{\ell^{2}(\hbar\mathbb{Z}^{n})}^{2}+\left\|\hbar^{-\alpha}(-\mathcal{L}_{\hbar})^{\frac{\alpha}{2}} u(t, \cdot)\right\|_{\ell^{2}(\hbar\mathbb{Z}^{n})}^{2}+\left\|m^{\frac{1}{2}}(\cdot) u(t, \cdot)\right\|_{\ell^{2}(\hbar\mathbb{Z}^{n})}^{2}.
	\end{equation}
	Then from \eqref{wellm}, it follows that
		\begin{eqnarray}
		\partial_{t} \sqrt{E(t)}&=&\frac{1}{2\sqrt{E(t)}}\partial_{t}E(t)\nonumber\\
		&=& \frac{1}{\sqrt{E(t)}}\operatorname{Re}\left(  f(t, \cdot), \partial_{t} u(t, \cdot)\right)_{\ell^{2}(\hbar\mathbb{Z}^{n})}\nonumber\\
		&\leq&\frac{1}{\sqrt{E(t)}}\left|\left(  f(t, \cdot), \partial_{t} u(t, \cdot)\right)_{\ell^{2}(\hbar\mathbb{Z}^{n})}\right|\nonumber\\
		&\leq&\frac{1}{\sqrt{E(t)}}\|f(t,\cdot)\|_{\ell^{2}(\hbar\mathbb{Z}^{n})}\|\partial_{t} u(t, \cdot)\|_{\ell^{2}(\hbar\mathbb{Z}^{n})}\nonumber\\
		&\leq&\|f(t,\cdot)\|_{\ell^{2}(\hbar\mathbb{Z}^{n})},
	\end{eqnarray}
	and so upon integrating, we obtain
	\begin{equation}
		\sqrt{E(t)}\leq \sqrt{E(0)}+\int_{0}^{t}\left\|f(s,\cdot)\right\|_{\ell^{2}(\hbar\mathbb{Z}^{n})}\mathrm{d}s.
	\end{equation}
Applying the Young's inequality followed by the H\"{o}lder's inequality, we get
\begin{eqnarray}\label{enrgyest}
	E(t)&\leq& E(0)+\left(\int_{0}^{t}\left\|f(s,\cdot)\right\|_{\ell^{2}(\hbar\mathbb{Z}^{n})}\mathrm{d}s\right)^{2}\nonumber\\
&\lesssim&	E(0)+\int_{0}^{t}\left\|f(s,\cdot)\right\|^{2}_{\ell^{2}(\hbar\mathbb{Z}^{n})}\mathrm{d}s,
\end{eqnarray}
for all $t\in[0,T]$, with constant depending on $T$ but not on $\hbar$.
Using Definition \ref{fraclap} and the Parseval's identity, we get
\begin{eqnarray}\label{bdd}
	\left\|\hbar^{-\alpha}(-\mathcal{L}_{\hbar})^{\frac{\alpha}{2}} u_{0}( \cdot)\right\|_{\ell^{2}(\hbar\mathbb{Z}^{n})}^{2}&=&(\hbar^{-\alpha}(-\mathcal{L}_{\hbar})^{\frac{\alpha}{2}} u_{0}( \cdot),\hbar^{-\alpha}(-\mathcal{L}_{\hbar})^{\frac{\alpha}{2}} u_{0}( \cdot))_{\ell^{2}(\hbar\mathbb{Z}^{n})}\nonumber\\
	&=&\hbar^{-2\alpha}\int\limits_{\mathbb{T}_{\hbar}^{n}}\left|2\sin(\pi\hbar\theta)\right|^{2\alpha}|\widehat{u_{0}}(\theta)|^{2}\mathrm{d} \theta\nonumber\\
	&=&\hbar^{-2\alpha}\int\limits_{\mathbb{T}_{\hbar}^{n}}\left[\sum\limits_{i=1}^{n}4\sin^{2}\left(\pi\hbar\theta_{i}\right)\right]^{\alpha}|\widehat{u_{0}}(\theta)|^{2}\mathrm{d} \theta\nonumber\\
	&\leq&\hbar^{-2\alpha}(4n)^{\alpha}\int\limits_{\mathbb{T}_{\hbar}^{n}}|\widehat{u_{0}}(\theta)|^{2}\mathrm{d} \theta\nonumber\\
	&\lesssim&\hbar^{-2\alpha}\left\| u_{0}\right\|_{\ell^{2}(\hbar\mathbb{Z}^{n})}^{2}.
\end{eqnarray}
Furthermore, the hypothesis that $m\in \ell^{\infty}(\hbar\mathbb{Z}^{n})$ gives
	\begin{equation}\label{mest}
		\left\|m^{\frac{1}{2}}(\cdot) u_{0}( \cdot)\right\|_{\ell^{2}(\hbar\mathbb{Z}^{n})}^{2}\lesssim \|m\|_{\ell^{\infty}(\hbar\mathbb{Z}^{n})}\left\|u_{0}\right\|_{\ell^{2}(\hbar\mathbb{Z}^{n})}^{2}.
	\end{equation}
	Combining the inequalities  \eqref{enrgyest}, \eqref{bdd}, and \eqref{mest}  with equality \eqref{enrgy}, we get 
	\begin{multline}\label{utest}
		\left\|\partial_{t} u(t, \cdot)\right\|_{\ell^{2}(\hbar\mathbb{Z}^{n})}^{2} \\
		\lesssim\left(\left\|u_{1}\right\|_{\ell^{2}(\hbar\mathbb{Z}^{n})}^{2}+\hbar^{-2\alpha}\left\| u_{0}\right\|_{\ell^{2}(\hbar\mathbb{Z}^{n})}^{2}+\|m\|_{\ell^{\infty}(\hbar\mathbb{Z}^{n})}\left\|u_{0}\right\|_{\ell^{2}(\hbar\mathbb{Z}^{n})}^{2}
		+\|f\|^{2}_{L^{2}([0,T];\ell^{2}(\hbar\mathbb{Z}^{n}))}\right),
	\end{multline}
	and
	\begin{multline}\label{normm}
		\left\|m^{\frac{1}{2}}(\cdot) u(t, \cdot)\right\|_{\ell^{2}(\hbar\mathbb{Z}^{n})}^{2}\\ \lesssim\left(\left\|u_{1}\right\|_{\ell^{2}(\hbar\mathbb{Z}^{n})}^{2}+\hbar^{-2\alpha} \left\| u_{0}\right\|_{\ell^{2}(\hbar\mathbb{Z}^{n})}^{2}+\|m\|_{\ell^{\infty}(\hbar\mathbb{Z}^{n})}\left\|u_{0}\right\|_{\ell^{2}(\hbar\mathbb{Z}^{n})}^{2}
		+\|f\|^{2}_{L^{2}([0,T];\ell^{2}(\hbar\mathbb{Z}^{n}))}\right),
	\end{multline}
	for all $t\in(0,T]$, with constants independent of $\hbar$.
	\medskip 
	
	Observe that to prove  estimate \eqref{uestt}, it remains to estimate  $\left\| u(t, \cdot)\right\|_{\ell^{2}(\hbar\mathbb{Z}^{n})}^{2}$.  Applying the Fourier transform to \eqref{mainpde} with respect to $k\in\hbar\mathbb{Z}^{n}$, we get
	\begin{equation}\label{ftpde}
		\left\{\begin{array}{l}
			\widehat{u}_{tt} (t, \theta)+\beta^{2}\widehat{u}(t,\theta)=\widehat{g}(t,\theta)+\widehat{f}(t,\theta) , \quad \text { with } (t,\theta) \in(0, T]\times \mathbb{T}_{\hbar}^{n}, \\
			\widehat{u}(0, \theta)=\widehat{u}_{0}(\theta), \quad \theta \in  \mathbb{T}_{\hbar}^{n}, \\
			\widehat{u}_{t}(0, \theta)=\widehat{u}_{1}(\theta), \quad \theta \in  \mathbb{T}_{\hbar}^{n},
		\end{array}\right.
	\end{equation}
	where $\beta^{2}=\hbar^{-2\alpha}|2\sin(\pi\hbar\theta)|^{2\alpha}$ and  $g(t,k)=-m(k)u(t,k)$. 
 If we denote
\begin{equation}
	v(t):=\widehat{u}(t,\theta),\quad  f(t):=\widehat{f}(t,\theta),\quad g(t):=\widehat{g}(t,\theta) ,
\end{equation}
and 
\begin{equation}
	v_{0}:=\widehat{u}_{0}(\theta),\quad v_{1}:=\widehat{u}_{1}(\theta),
\end{equation}
then  \eqref{ftpde} becomes
\begin{equation}\label{ftpde1}
	\left\{\begin{array}{l}
		v^{\prime\prime}(t)+\beta^{2}v(t)=g(t)+f(t), \quad \text { with } t \in(0, T], \\
		v(0)=v_{0}, \quad v^{\prime}(0)=v_{1},
	\end{array}\right.
\end{equation}
with $\beta\geq0$  . By solving the homogeneous part of \eqref{ftpde1}, and using the Duhamel's principle (see, e.g. \cite{evans}), we get
\begin{equation}\label{vsol}
	v(t)=\cos \left(\beta t\right) v_{0}+\frac{\sin \left(\beta t\right)}{\beta } v_{1} 
	+\int_{0}^{t} \frac{\sin \left(\beta(t-s)\right)}{\beta}\left(g(s)+f(s)\right)  \mathrm{d} s.
\end{equation}
Using the following inequalities:
\begin{equation}
	|\cos(\beta t)|\leq 1 \text{~for all~} t\in[0,1]\text{~and~} |\sin(\beta t)|\leq 1,
\end{equation}
for large values of $\beta t$, while $|\sin(\beta t)|\leq \beta t\leq \beta T$ for small values of $\beta t, $ the equality \eqref{vsol} gives
\begin{equation}
	|v(t)|\lesssim|v_{0}|+|v_{1}|+\int\limits_{0}^{t}|g(s)|ds+\int\limits_{0}^{t}|f(s)|ds.
\end{equation}
	Now taking the $L^{2}$-norm, we get
	\begin{multline}
		\|\widehat{u}(t, \cdot)\|^{2}_{L^{2}(\mathbb{T}_{\hbar}^{n})} \lesssim\left\|\widehat{u}_{0}\right\|^{2}_{L^{2}(\mathbb{T}_{\hbar}^{n})}+\left\|\widehat{u}_{1}\right\|^{2}_{L^{2}(\mathbb{T}_{\hbar}^{n})}+\int_{0}^{t}\|\widehat{g}(s, \cdot)\|^{2}_{L^{2}(\mathbb{T}_{\hbar}^{n})} \mathrm{d} s\\+\int_{0}^{t}\|\widehat{f}(s, \cdot)\|^{2}_{L^{2}(\mathbb{T}_{\hbar}^{n})} \mathrm{d} s.
	\end{multline}
	Using the  Plancherel formula, we obtain
	\begin{equation}\label{normu}
		\|u(t, \cdot)\|^{2}_{\ell^{2}(\hbar\mathbb{Z}^{n})} \lesssim\left\|u_{0}\right\|^{2}_{\ell^{2}(\hbar\mathbb{Z}^{n})}+\left\|u_{1}\right\|^{2}_{\ell^{2}(\hbar\mathbb{Z}^{n})}+\int_{0}^{t}\|g(s, \cdot)\|^{2}_{\ell^{2}(\hbar\mathbb{Z}^{n})} \mathrm{d} s+\int_{0}^{t}\|f(s, \cdot)\|^{2}_{\ell^{2}(\hbar\mathbb{Z}^{n})} \mathrm{d} s,
	\end{equation}
	for all $t\in(0,T]$, where
	\begin{eqnarray}\label{normf}
		\|g(t, \cdot)\|^{2}_{\ell^{2}(\hbar\mathbb{Z}^{n})}=\|m(\cdot)u(t, \cdot)\|^{2}_{\ell^{2}(\hbar\mathbb{Z}^{n})}\lesssim\left\|m\right\|_{\ell^{\infty}(\hbar\mathbb{Z}^{n})} \left\|m^{\frac{1}{2}}(\cdot)u(t, \cdot)\right\|^{2}_{\ell^{2}(\hbar\mathbb{Z}^{n})}.	
	\end{eqnarray}
	Combining the inequalities \eqref{normm}, \eqref{normu} and \eqref{normf}, we get
	\begin{multline}\label{uest}
		\|u(t, \cdot)\|^{2}_{\ell^{2}(\hbar\mathbb{Z}^{n})} \lesssim\left\|u_{0}\right\|^{2}_{\ell^{2}(\hbar\mathbb{Z}^{n})}+\left\|u_{1}\right\|^{2}_{\ell^{2}(\hbar\mathbb{Z}^{n})}+T\left\|m\right\|_{\ell^{\infty}(\hbar\mathbb{Z}^{n})}\left(\left\|u_{1}\right\|_{\ell^{2}(\hbar\mathbb{Z}^{n})}^{2}+\right.\\\left.\hbar^{-2\alpha} \left\| u_{0}\right\|_{\ell^{2}(\hbar\mathbb{Z}^{n})}^{2}+\|m\|_{\ell^{\infty}(\hbar\mathbb{Z}^{n})}\left\|u_{0}\right\|_{\ell^{2}(\hbar\mathbb{Z}^{n})}^{2}+\|f\|^{2}_{L^{2}([0,T];\ell^{2}(\hbar\mathbb{Z}^{n}))}\right)\\
		+\|f\|^{2}_{L^{2}([0,T];\ell^{2}(\hbar\mathbb{Z}^{n}))},
	\end{multline}
	for all $t\in(0,T]$. Then from \eqref{utest} and \eqref{uest}, it follows that
	\begin{multline}
		\|u(t, \cdot)\|^{2}_{\ell^{2}(\hbar\mathbb{Z}^{n})}+\left\|\partial_{t} u(t, \cdot)\right\|_{\ell^{2}(\hbar\mathbb{Z}^{n})}^{2}\lesssim\left(1+\|m\|_{\ell^{\infty}(\hbar\mathbb{Z}^{n})}\right)\times\\\left[\left(\hbar^{-2\alpha}+\|m\|_{\ell^{\infty}(\hbar\mathbb{Z}^{n})}\right)\left\|u_{0}\right\|_{\ell^{2}(\hbar\mathbb{Z}^{n})}^{2}+\left\|u_{1}\right\|_{\ell^{2}(\hbar\mathbb{Z}^{n})}^{2}+\|f\|^{2}_{L^{2}([0,T];\ell^{2}(\hbar\mathbb{Z}^{n}))}\right],
	\end{multline}
	for all $t\in(0,T]$, with constant depending on $T$ but not on $\hbar$. This gives the  required estimate \eqref{uestt} from which the
 uniqueness  follows immediately.
\end{proof}
\section{Limit $\hbar\to 0$}\label{limit}
\begin{proof}[Proof of Theorem \ref{cgt}]
Consider two Cauchy problems:
\begin{equation}\label{CP1}
	\left\{\begin{array}{l}
		\partial_{t}^{2} u(t, k)+\hbar^{-2\alpha}(-\mathcal{L}_{\hbar})^{\alpha} u(t, k)+m(k) u(t, k)=0, \quad  (t,k) \in(0,T]\times\hbar\mathbb{Z}^{n}, \\
		u(0, k)=u_{0}(k), \quad k\in\hbar\mathbb{Z}^{n},\\
		\partial_{t} u(0, k)=u_{1}(k),\quad k\in\hbar\mathbb{Z}^{n},
	\end{array}\right.
\end{equation}
and
\begin{equation}\label{CP2}
	\left\{
	\begin{array}{ll}
		\partial^{2}_{t}v(t,x)+(-\mathcal{L})^{\alpha}v(t,x)+m(x)v(t,x)=0, \quad (t,x) \in(0,T]\times\mathbb{R}^{n},\\
		v(0,x)=u_{0}(x),\quad x\in\mathbb{R}^{n},\\
		\partial_{t}v(0,x)=u_1(x),\quad x\in\mathbb{R}^{n},
	\end{array}
	\right.
\end{equation}
where $\left(-\mathcal{L}\right)^{\alpha}$ is the usual fractional Laplacian on $\mathbb{R}^{n}$. Here the initial data of the Cauchy problem \eqref{CP1} is the evaluation of the initial data from \eqref{CP2} on the lattice $\hbar\mathbb{Z}^{n}$. From the equations \eqref{CP1} and \eqref{CP2}, 
denoting $w:=u-v$,  we get
\begin{equation}\label{CPF}
	\left\{
	\begin{array}{ll}
		\partial^{2}_{t}w(t,k)+\hbar^{-2\alpha}(-\mathcal{L}_{\hbar})^{\alpha}w(t,k)+m(k)w(t,k)=\left((-\mathcal{L})^{\alpha}-\hbar^{-2\alpha}(-\mathcal{L}_{\hbar})^{\alpha}\right)v(t,k),\\
		w(0,k)=0,\quad k\in\hbar\mathbb{Z}^{n},\\
		\partial_{t}w(0,k)=0,\quad k\in\hbar\mathbb{Z}^{n}.               
	\end{array}
	\right.
\end{equation}
Since $w_{0}=w_{1}=0$, applying the Theorem \ref{wellpo} for the above Cauchy problem and using the estimate \eqref{uestt}, we get   
\begin{multline}\label{cgt1}
\left\|w(t)\right\|^{2}_{\ell^{2}(\hbar\mathbb{\mathbb{Z}}^{n})}+\left\|\partial_{t}w(t)\right\|^{2}_{\ell^{2}(\hbar\mathbb{\mathbb{Z}}^{n})}\lesssim\\ (1+\left\|m\right\|_{\ell^{\infty}(\hbar\mathbb{Z}^{n})})\left\|\left((-\mathcal{L})^{\alpha}-\hbar^{-2\alpha}(-\mathcal{L}_{\hbar})^{\alpha}\right)v(t)\right\|^{2}_{\ell^{2}(\hbar\mathbb{Z}^{n})},
\end{multline}
for all $t\in[0,T].$ Now we will estimate the term $\left\|\left((-\mathcal{L})^{\alpha}-\hbar^{-2\alpha}(-\mathcal{L}_{\hbar})^{\alpha}\right)v(t)\right\|_{\ell^{2}(\hbar\mathbb{Z}^{n})}.$ 
Using the Plancherel formula, we have
\begin{equation}\label{fnm}
	\left\|\left((-\mathcal{L})^{\alpha}-\hbar^{-2\alpha}(-\mathcal{L}_{\hbar})^{\alpha}\right)v(t)\right\|_{\ell^{2}(\hbar\mathbb{Z}^{n})}=\left\|(|2\pi\theta|^{2\alpha}-\hbar^{-2\alpha}|2\sin(\pi\hbar\theta)|^{2\alpha})\widehat{v}(t)\right\|_{L^{2}(\mathbb{T}_{\hbar}^{n})}.
\end{equation}
Since  $|\cdot|^{\alpha}:\mathbb{R}\to \mathbb{R}$ is $\alpha$-H\"older continuous for $0<\alpha<1$,  we have the following inequality
\begin{equation}
	\left||x|^{2\alpha}-|y|^{2\alpha}\right|\lesssim \left||x|^{2}-|y|^{2}\right|^{\alpha}, \quad x,y\in\mathbb{R}^{n},
\end{equation} 
 and  $|x|=\sqrt{x_{1}^{2}+\dots+x_{n}^{2}}$. Now, using the above inequality and the Taylor expansion for $\sin^{2}(\pi\hbar\theta_{j})$, we get
\begin{eqnarray}
	\left||2\pi\theta|^{2\alpha}-\hbar^{-2\alpha}|2\sin(\pi\hbar\theta)|^{2\alpha}\right|
	&\lesssim&
	\left|\sum\limits_{j=1}^{n}4\pi^2\theta^{2}_{j}-\hbar^{-2}\sum\limits_{j=1}^{n}4\sin^{2}(\pi\hbar\theta_{j})\right|^{\alpha}\nonumber\\
	&=&\left|\sum\limits_{j=1}^{n}\left[4\pi^2\theta^{2}_{j}-\hbar^{-2}4\left(\pi^{2}\hbar^{2}\theta_{j}^{2}-\frac{\pi^{4}}{3}\hbar^{4}\theta_{j}^{4}\cos\left(2\xi_{j}\right)\right)\right]\right|^{\alpha}\nonumber\\
	&=&\left(\frac{4\pi^{4}\hbar^{2}}{3}\right)^{\alpha}\left|\sum\limits_{j=1}^{n}\theta_{j}^{4}\cos(2\xi_{j})\right|^{\alpha}\nonumber\\
	&\lesssim&\hbar^{2\alpha}\left[\sum\limits_{j=1}^{n}\theta^{4}_{j}\right]^{\alpha}\nonumber\\
	&\lesssim&\hbar^{2\alpha}|\theta|^{4\alpha},\quad \theta\in\mathbb{T}_{\hbar}^{n},
\end{eqnarray}
where $|\theta|^{4\alpha}=\left[\sum\limits_{j=1}^{n}\theta^{2}_{j}\right]^{2\alpha}$ and $\xi_{j}\in(0,\pi\hbar\theta_{j})$ or $(\pi\hbar\theta_{j},0)$ depending on the sign of $\theta_{j}$.  Since $u_{0}\in H^{4\alpha}(\mathbb{R}^{n})$, we have $v\in H^{4\alpha}(\mathbb{R}^{n})$ by Theorem \ref{fracklien}. Now, combining the above estimate with  \eqref{fnm}, we get
\begin{equation}\label{fnmm}
\left\|\left((-\mathcal{L})^{\alpha}-\hbar^{-2\alpha}(-\mathcal{L}_{\hbar})^{\alpha}\right)v(t)\right\|_{\ell^{2}(\hbar\mathbb{Z}^{n})}\lesssim \hbar^{2\alpha}\||\theta|^{4\alpha}\widehat{v}(t)\|_{L^{2}(\mathbb{T}_{\hbar}^{n})} ,
\end{equation}
Using \eqref{cgt1} and \eqref{fnmm}, we get $\|w(t)\|_{\ell^{2}\left(\hbar \mathbb{Z}^{n}\right)}^{2}+\left\|\partial_{t} w(t)\right\|_{\ell^{2}\left(\hbar \mathbb{Z}^{n}\right)}^{2} \rightarrow 0$ as $\hbar \rightarrow 0$. Hence $\|w(t)\|_{\ell^{2}\left(\hbar \mathbb{Z}^{n}\right)} \rightarrow 0$ and $\left\|\partial_{t} w(t)\right\|_{\ell^{2}\left(\hbar \mathbb{Z}^{n}\right)} \rightarrow 0$ as $\hbar \rightarrow 0$. This finishes the proof of Theorem \ref{cgt}.
\end{proof}

\bibliographystyle{abbrv}
\bibliography{discrete_fractional_laplace}

\begin{thebibliography}{10}

\bibitem{abuteen}
E.~Abuteen, A.~Freihat, M.~Al-Smadi, H.~Khalil, and R.~A. Khan.
\newblock Approximate series solution of nonlinear, fractional klein-gordon
  equations using fractional reduced differential transform method.
\newblock {\em Journal of Mathematics and Statistics}, 12(1):23--33, Mar. 2016.

\bibitem{fracklien}
A.~Altybay, M.~Ruzhansky, M.~E. Sebih, and N.~Tokmagambetov.
\newblock Fractional {K}lein-{G}ordon equation with singular mass.
\newblock {\em Chaos Solitons Fractals}, 143:Paper No. 110579, 6, 2021.

\bibitem{fracschro}
A.~Altybay, M.~Ruzhansky, M.~E. Sebih, and N.~Tokmagambetov.
\newblock Fractional {S}chr\"{o}dinger equation with singular potentials of
  higher order.
\newblock {\em Rep. Math. Phys.}, 87(1):129--144, 2021.

\bibitem{altu}
A.~Arda, R.~Sever, and C.~Tezcan.
\newblock Analytical solutions to the klein{\textendash}gordon equation with
  position-dependent mass for q-parameter pöschl{\textendash}teller potential.
\newblock {\em Chinese Physics Letters}, 27(1):010306, jan 2010.

\bibitem{kibti}
L.~N.~A. Botchway, P.~Ga\"{e}l~Kibiti, and M.~Ruzhansky.
\newblock Difference equations and pseudo-differential operators on {$\Bbb
  Z^n$}.
\newblock {\em J. Funct. Anal.}, 278(11):108473, 41, 2020.

\bibitem{mari}
M.~Chatzakou, M.~Ruzhansky, and N.~Tokmagambetov.
\newblock Fractional {K}lein-{G}ordon equation with singular mass. {II}:
  {H}ypoelliptic case.
\newblock {\em Complex Var. Elliptic Equ.}, 67(3):615--632, 2022.

\bibitem{roncal}
O.~Ciaurri, L.~Roncal, P.~R. Stinga, J.~L. Torrea, and J.~L. Varona.
\newblock Nonlocal discrete diffusion equations and the fractional discrete
  {L}aplacian, regularity and applications.
\newblock {\em Adv. Math.}, 330:688--738, 2018.

\bibitem{dasgupta2021discrete}
A.~Dasgupta, M.~Ruzhansky, and A.~Tushir.
\newblock Discrete time-dependent wave equations {I}. {S}emiclassical analysis.
\newblock {\em Journal of Differential Equations}, 317:89--120, 2022.

\bibitem{desouza}
A.~de~Souza~Dutra and C.-S. Jia.
\newblock Classes of exact {K}lein-{G}ordon equations with spatially dependent
  masses: regularizing the one-dimensional inversely linear potential.
\newblock {\em Phys. Lett. A}, 352(6):484--487, 2006.

\bibitem{ege}
S.~M. Ege and E.~Misirli.
\newblock Solutions of the space-time fractional foam drainage equation and the
  fractional klein-gordon equation by use of modified kudryashov method.
\newblock 2014.

\bibitem{evans}
L.~C. Evans.
\newblock {\em Partial differential equations}, volume~19 of {\em Graduate
  Studies in Mathematics}.
\newblock American Mathematical Society, Providence, RI, 1998.

\bibitem{fracthought}
N.~Garofalo.
\newblock Fractional thoughts.
\newblock In {\em New developments in the analysis of nonlocal operators},
  volume 723 of {\em Contemp. Math.}, pages 1--135. Amer. Math. Soc.,
  [Providence], RI, [2019] \copyright 2019.

\bibitem{garra}
R.~Garra, E.~Orsingher, and F.~Polito.
\newblock Fractional klein-gordon equation for linear dispersive phenomena:
  Analytical methods and applications.
\newblock pages 1--6, 2014.

\bibitem{gepr}
K.~A. Gepreel and M.~S. Mohamed.
\newblock Analytical approximate solution for nonlinear space{\textemdash}time
  fractional klein{\textemdash}gordon equation.
\newblock {\em Chinese Physics B}, 22(1):010201, jan 2013.

\bibitem{ghosh}
U.~Ghosh, J.~Banerjee, S.~Sarkar, and S.~Das.
\newblock Fractional klein--gordon equation composed of jumarie fractional
  derivative and its interpretation by a smoothness parameter.
\newblock {\em Pramana}, 90(6):1--10, 2018.

\bibitem{gol}
A.~K. Golmankhaneh, A.~K. Golmankhaneh, and D.~Baleanu.
\newblock On nonlinear fractional klein--gordon equation.
\newblock {\em Signal Processing}, 91(3):446--451, 2011.

\bibitem{hua}
Y.~Huang and A.~Oberman.
\newblock Numerical methods for the fractional {L}aplacian: a finite
  difference--quadrature approach.
\newblock {\em SIAM J. Numer. Anal.}, 52(6):3056--3084, 2014.

\bibitem{kirk}
K.~Kirkpatrick, E.~Lenzmann, and G.~Staffilani.
\newblock On the continuum limit for discrete {NLS} with long-range lattice
  interactions.
\newblock {\em Comm. Math. Phys.}, 317(3):563--591, 2013.

\bibitem{markus}
M.~Klein and E.~Rosenberger.
\newblock Harmonic approximation of difference operators.
\newblock {\em J. Funct. Anal.}, 257(11):3409--3453, 2009.

\bibitem{ten}
M.~Kwa\'{s}nicki.
\newblock Ten equivalent definitions of the fractional {L}aplace operator.
\newblock {\em Fract. Calc. Appl. Anal.}, 20(1):7--51, 2017.

\bibitem{duarte}
M.~D. Ortigueira.
\newblock {\em Fractional calculus for scientists and engineers}, volume~84 of
  {\em Lecture Notes in Electrical Engineering}.
\newblock Springer, Dordrecht, 2011.

\bibitem{Rab10}
V.~Rabinovich.
\newblock Exponential estimates of solutions of pseudodifferential equations on
  the lattice {$(h\Bbb Z)^n$}: applications to the lattice {S}chr\"odinger and
  {D}irac operators.
\newblock {\em J. Pseudo-Differ. Oper. Appl.}, 1(2):233--253, 2010.

\bibitem{Rab13}
V.~Rabinovich.
\newblock Wiener algebra of operators on the lattice {$(\mu \Bbb{Z})^n$}
  depending on the small parameter {$\mu>0$}.
\newblock {\em Complex Var. Elliptic Equ.}, 58(6):751--766, 2013.

\bibitem{Rab18}
V.~S. Rabinovich and S.~Roch.
\newblock Pseudodifference operators on weighted spaces, and applications to
  discrete {S}chr\"odinger operators.
\newblock {\em Acta Appl. Math.}, 84(1):55--96, 2004.

\bibitem{RT-Birk}
M.~Ruzhansky and V.~Turunen.
\newblock On the {F}ourier analysis of operators on the torus.
\newblock In {\em Modern trends in pseudo-differential operators}, volume 172
  of {\em Oper. Theory Adv. Appl.}, pages 87--105. Birkh\"auser, Basel, 2007.

\bibitem{RT-JFAA}
M.~Ruzhansky and V.~Turunen.
\newblock Quantization of pseudo-differential operators on the torus.
\newblock {\em J. Fourier Anal. Appl.}, 16(6):943--982, 2010.

\bibitem{samko1993fractional}
S.~G. Samko, A.~A. Kilbas, O.~I. Marichev, et~al.
\newblock {\em Fractional integrals and derivatives}, volume~1.
\newblock Gordon and breach science publishers, Yverdon Yverdon-les-Bains,
  Switzerland, 1993.

\bibitem{wang1}
B.-Q. Wang, Z.-W. Long, C.-Y. Long, and S.-R. Wu.
\newblock Klein-{G}ordon oscillator with position-dependent mass in the
  rotating cosmic string spacetime.
\newblock {\em Modern Phys. Lett. A}, 33(4):1850025, 11, 2018.

\bibitem{wang}
Z.~Wang, Z.~Long, C.~Long, and L.~Wang.
\newblock Analytical solutions of position-dependent mass klein--gordon
  equation for unequal scalar and vector yukawa potentials.
\newblock {\em Indian Journal of Physics}, 89(10):1059--1064, 2015.

\bibitem{culha}
S.~Çulha and A.~Daşcıoğlu.
\newblock Analytic solutions of the space–time conformable fractional
  klein–gordon equation in general form.
\newblock {\em Waves in Random and Complex Media}, 29(4):775--790, 2019.

\end{thebibliography}

\end{document}